\documentstyle[dina4,mssymb]{article}

\newenvironment{Proof}{\noindent {\sc Proof:}\begin{em}}%
{\end{em}\hfill \rule{2mm}{2mm}}

\newtheorem{twierdz}{Theorem }[section]
\newtheorem{stwierdz}[twierdz]{Proposition }

\newtheorem{definicja}{Definition }[section]

\newtheorem{uwaga}[twierdz]{Remark} 
\newtheorem{hipoteza}{Conjecture}

\def\ot{\otimes}
\def\ra{\longmapsto}

\def\cch{\mbox{\gt H}}

\def\cb{{\cal B}}

\def\cd{{\cal D}}
\def\ce{{\cal E}}

\def\ch{{\cal H}}

\def\cw{{\cal W}}

\def\idt{{\bf 1}}

\def\1n{^{(1)}}
\def\2n{^{(2)}}
\def\3n{^{(3)}}

\def\t{\tilde}

\def\ba{\begin{array}}
\def\ea{\end{array}}  
\def\be{\begin{equation}}
\def\ee{\end{equation}}
\def\bdm{\begin{displaymath}}
\def\edm{\end{displaymath}}
\def\beqnar{\begin{eqnarray}}
\def\eeqnar{\end{eqnarray}}
\def\beqna*{\begin{eqnarray*}}
\def\eeqna*{\end{eqnarray*}}

\def\nsl{/}

\def\ch{{\cal H}}
\def\<{\langle}
\def\>{\rangle}
\newfont{\gt}{eufb10 scaled 1000}

\begin{document}
\title{Wick algebras approach to physics of 2D systems
\footnote{Submited to the Proceedings of the Fifth International School 
on Theoretical Physics "Symmetry and Structural Properties of Condensed Matter",
1998, Zajaczkowo, Poland, World Scientific P.C.}}

\author{{\sc Roman Gielerak}\thanks{Technical University of Zielona G{\.o}ra, 
Institute of T. Physics, Poland}\\
		University of Wroc\l aw\\
		Institute of Theoretical Physics\\
		50-204 Wroc\l aw, Poland\\
		{\sc Robert Ra\l owski}\\
		Wroc\l aw University of Technology\\
		Institute of Mathematics\\
		Wroc\l aw, Poland}

\maketitle
\begin{abstract}
The notion of anyonic Wick algebras is introduced and the 
corresponding second quantisation procedure is discussed. Applications to
the physics of 2D anyonic matter are presented. In particular the existence of 
thermodynamics in the infinite volume for the case of Leinaas-Myrheim 
potentials is presented. Additionally a new class of states corresponding 
to quasifree thermal states on some Wick algebras is described.
\end{abstract}

\section{Introduction}
It is well known that after lowering the dimension of the space to two, the
corresponding notion of statistics of the system of quantum particles
admits a significant generalisation~\cite{Green,Lerda,Leinaas,W1,W2}.

The corresponding wave functions $\Psi_n(x_1,...,x_n)$ form a Hilbert space
where unitary representation of the so called braid group $B_n$ acts (where $n$ 
is
equal to the number of particles). The braid group $B_n$ known as Artin group 
is generated by $\{1,t_1,...,t_{n-1}\}$ fulfilling:
$$
t_it_j=t_jt_i\mbox{ for }|i-j|\ge 2\mbox{ and }t_{i+1}t_it_{i+1}=t_it_{i+1}t_i
\mbox{ for }i=1,2,...,n-2
$$
In contrast to the familiar permutation group $S_n$ the group $B_n$ (for $n>1$) 
is
infinite. The braid group symmetry of the wave functions under exchange of two
particles leads to the idea of describing such systems of particles in terms of 
the
cyclic representations of the underlying algebra of the exotic commutation 
relations.
Such algebras are called Wick algebras in the 
literature~\cite{BS1,BS2,Green,LM,JSW}.

In the present exposition we shall focus our attention on a particular class
of Wick algebras in which the corresponding commutation relations lead to the
realisation of the groups $B_n$ in the corresponding Fock modules as the group
$S_n$. Such Wick algebras are called anyonic type Wick algebras.
A special interest in this class of Wick algebras comes from the fact that
systems of particles interacting throughout the well known
Leinaas-Myrheim hamiltonians~\cite{Lerda,Leinaas} can be described by such 
algebras.
Additionally,
the arguments outlined below provide us with some rigorous results penetrating
the physics of 2D anyonic systems. For a complete exposition of the proofs 
we refer to~\cite{GR}.

\section{Anyonic type Wick algebras and the corresponding Fock modules}
Let $\cch$ be complex separable Hilbert space and let $T(\cch,\cch^*)$ be the
corresponding full tensor algebra over the space $\cch\oplus\cch^*$ where 
$\cch^*$
stands for the conjugate (via the Riesz representation Theorem).
Let $\t{T}\in\cb(\cch\ot\cch^*,\cch^*\ot\cch)$ be given and let $J_T$ be an 
ideal 
generated by $f^*\ot g-\t{T}(g\ot f^*)$. Then the algebra over the space $\cch$ 
and
generated by the commutation relations given by $T$, where $T$ is 
a partial conjugate of $\t{T}$ is called a Wick algebra generated by the
commutation relations $T$. In other words the corresponding Wick algebra is
equal to the quotient $T(\cch,\cch^*)\nsl J_{\t{T}}$
We denote the obtained Wick algebra as $\cw(\cch,T)$.

A Fock module $\Gamma_T^0(\cch)$ for given $\cw(\cch,T)$ is a cyclic left module
with cyclic vector $\Omega$ obeing $a(f)\Omega=0$ for any $f^*\in\cch^*$, where 
we
have denoted by $a(f)$ the action of the tensor $f\in\cw(\cch,T)$ 
in $\Gamma_T(\cch)$. From the definition of the Fock module $\Gamma_T^0(\cch)$
it follows that there exists an unique sesquilinear form $\<\cdot,\cdot\>_T$
on 
$\Gamma_T^0(\cch)$ such that the adjoint $a^+(f)$ of $a(f)$ is equal to 
$a(f^*)$,
i.e. $\<x,a(f)y\>_T=\<a(f)^+x,y\>_T$ for suitable $x,y\in\Gamma_T(\cch)$.
By using the corresponding commutation relations and the property $a(f)\Omega=0$
it follows that
\beqna*
\<a^+(f_1)..a^+(f_n)\Omega,a^+(g_1)..a^+(g_m)\Omega\>_T
&=&\delta_{n,m}
\<f_1\ot..\ot f_n,P_n(T)g_1\ot..\ot g_m\>_T\\ 
\eeqna*
where $P_n(T)$ is defined recursively as $P_0(T)=\idt$, $P_1(T)=\idt$ and
$P_{n+1}(T)=(P_n(T)\ot\idt)(\idt+T_n+...T_n...T_1)$, where $T_i$ is an operator
$\idt\ot...\ot T\ot...\ot\idt$ acting in $\cch^{\ot n+1}$ for $i=1,..,n$.

Of particular interest to physics are situations in which $P_n(T)\ge0$.
The following sufficient conditions for this to hold are 
known~\cite{BS1,BS2,JSW}:
\begin{enumerate}
\item $T\ge0$
\item $\|T\|\le\frac{1}{2}$ and if $\|T\|<\frac{1}{2}$ then $kerP(T)=\{ 0\}$
\item $T$ fulfills YB relation $T_1T_2T_1=T_2T_1T_2$ in $\cch\ot\cch\ot\cch$
and $\|T\|\le1$. If moreover $\|T\|<1$ then $P_n(T)>0$ for all $n$.
\end{enumerate}
In the following we shall assume that one of the conditions $1-3$ holds
and then the metric completion of $\Gamma_T^0(\cch)$ denoted as $\Gamma_T(\cch)$
(after eventual division by the kernel of the $\bigoplus_{n\ge0}P_n$)
is called the Fock module of the algebra $\cw(\cch,T)$.

Let $h$ be positive selfadjoint operator acting in $\cch$ and such that
$Tr(e^{-\beta h})<\infty$ for any $0<\beta<\infty$. We call it kinetic
energy operator. According to the standard construction we will be
interested in lifting the semigroup $e^{-\beta h}$ to the corresponding 
Fock module $\Gamma_T(\cch)$.
\begin{stwierdz}
\begin{enumerate}
\item Let $A\in\cb(\cch)$ be selfadjoint. Define: 
$$
d_T\Gamma(A)=\bigoplus_{n\ge0}(\sum_{k=1}^n
\underbrace{\idt\ot...\ot A\ot...\ot \idt)}_{A\hbox{ is on the k-th place}}
$$
on $\bigoplus_{n\ge0}\cd(A)^{\ot n}$ then $d_T\Gamma(A)$ is symmetric iff 
$[T,d\Gamma_2(A)]=0$, where $d\Gamma_2(A)=\idt\ot A+A\ot\idt$.
Morever if $A$ is essentially self adjoint operator on $\cd(A)$ then 
$d_T\Gamma(A)$ is also essentially selfadjoint on 
$\bigoplus^{alg}_{n\ge0}\cd(A)^{\ot n}$.
\item Let $A\in\cb(\cch)$ be unitary. Define:
$\Gamma_T(A)=\bigoplus_{n\ge0}A^{\ot n}$ then $\Gamma_T(A)$ is unitary in
$\Gamma_T(\cch)$ iff $[T,d\Gamma_2(A)]=0$.
\end{enumerate}
\end{stwierdz}

For the standard systems like bosons or fermions, the macroscopic
thermodynamical properties of the corresponding quantum (non) interacting gases 
are encoded in the so called partition function. Let $h$ be the kinetic energy
operator and let $\Gamma_T(e^{-\beta h})$ be the corresponding density matrix
in the Fock module.

\begin{stwierdz}
Assume that $Tr(e^{-\beta h})<\infty$ and let $T$ be such that $P(T)>0$
(i.e. $ker P(T)=\{ 0\}$) then $\Gamma_T(e^{-\beta h})$ is of trace class on
$\Gamma_T(\cch)$ and:
$$
Tr_{\Gamma_T(\cch)}\Gamma_T(e^{-\beta 
h})=Tr_{\Gamma_{T=0}(\cch)}\Gamma_0(e^{-\beta h})
$$
\end{stwierdz}

Thus we see that the thermodynamics of the system of particles obeying 
$T$-commutation relations with $T$ such that $P(T)>0$ is the same as the  
thermodynamics of the so called system of "Maxwell-Boltzmann particles".
The nontrivial influence on the thermodynamic properties of the deformation
of standard commutation relations might appear only in the case when the  
corresponding metric form $P(T)$ has a nontrivial kernel.

For the so called $q$-commutation relations~\cite{BS1,BS2} this effect has been
detected in~\cite{Werner}. The strict positivity of $P(T)$ yields also that 
there is no nontrivial quadratic Wick ideals which is equivalent to the
statement that there is no (simple) commutation relation between the
corresponding creation and annihilation operators. This means in particular 
that there is no symmetry of the vector 
$\Psi(f_,...,f_n)=a^+(f_1)...a^+(f_n)\Omega$ under the exchange of $f_i\ra f_j$
and thus the system behaves like that in the case $T=0$.

In order to see how the deformation of standard commutation relations influence
the partition functions (and therefore the physics) we select a class
of Wick algebras for which (as we shall see in the next paragraph) this problem
seems to be tractable by standard methods of Quantum Statistical 
Mechanics~\cite{BR}.

\begin{definicja}
A Wick algebra $\cw(\cch,T)$ is called anyonic type Wick algebra iff the  
operator $T=T^*$ determining the commutation relations fulfills $YBE$ and
$(\idt-T)(\idt+\hat{T})=0$, where if $T=(t_{i,j}^{k,l})$ in some CONS $(e_n)$
in $\cch$ then $\hat{T}=(\hat{t}_{i,j}^{k,l})=t(_{j,l}^{i,k})$.
\end{definicja}

\begin{stwierdz}
Let $\cw(\cch,T)$ be an anyonic type Wick algebra. Then there exists a Fock
representation ($\Gamma_T(\cch,\Omega_T,a(f),a^+(f)$) of $\cw(\ch,T)$ such that 
the following
algebra of commutation relations is fulfilled:
$$
a(e_i)a^+(e_k)-\sum_{k,l}t_{i,j}^{k,l}a^+(e_k)a(e_l)=\delta_{i,j}\idt
$$
$$
a(e_i)a(e_j)-\sum_{k,l}t_{i,j}^{k,l}a(e_k)a(e_l)=0
$$
$$
a^+(e_i)a^+(e_j)-\sum_{k,l}t_{i,j}^{k,l}a^+(e_k)a^+(e_l)=0
$$
In particular it follows that $ker\<\cdot,\cdot\>_T=P(T)\Gamma_0(\cch)^{\perp}$,
where $\Gamma_0(\cch)$ is the free Fock module over $\cch$ and $P(T)$ is
an orthogonal projector.
\end{stwierdz}
\begin{Proof}
See~\cite{GR}
\end{Proof}

The following strategy for displaying the influence of the nonstardard 
commutation relations on the thermodynamical properties of the system under 
consideration will be realised in the next section.

Let $e^{-\beta h}$ be a trace class semigroup acting in $\cch$ and
let $\cw(\cch,T)$ be anyonic type Wick algebra. Denote by $\Gamma_T(\cch)$
the corresponding Fock module and let $\Gamma_{+\atop -}(\cch)$ denotes
the standard fermionic (-) or bosonic (+) Fock modules. Let $J^{+\atop-}$ be
some unitary isomorphism in between $\Gamma_T(\cch)$ and 
$\Gamma_{+\atop-}(\cch)$ (where we assumed that dim $\cch=\infty$).
Using the invariance of the trace under the unitary transformation and assuming
that 
$J^{+\atop-}\Gamma_T(e^{-\beta h})(J^{+\atop-})^{-1}=
exp(-\beta(d\Gamma_{+\atop-}(h)+R))$,
where $R$ is some tractable perturbation of the second quantised free
hamiltonian $d\Gamma_{+\atop-}(h)$, we arrive at:
$$
Z_T(h)\equiv Tr_{\Gamma_T(\cch)}\Gamma_T(e^{-\beta h})=
Tr_{\Gamma_{+\atop-}(\cch)}(e^{-\beta(d\Gamma_{+\atop-}(h)+R)})
$$

\section{Thermodynamics of the 2D Anyonic matter}
Let $\Lambda\subset R^d$ be a bounded region with $C^1$-piecewise boundary 
$\partial\Lambda$ and let for $\sigma\in C(\partial\Lambda)$, $\sigma(x)\ge0$,
$-\Delta^\sigma_\Lambda$ be a selfadjoint version of the Laplace operator
$-\Delta$ acting in the space $L^2(\Lambda,dx)$ and corresponding to the  
(classical) 
boundary condition $\sigma$, see~\cite{BR}. Let $\t{r}(x,y)=r(x-y)$ be such that
$r\in C^2(\Bbb{R}^d)$ and $\t{r}(x,y)+\t{r}(y,x)=2k\pi$, $k\in\Bbb{Z}$.
In the space $L^2(\Lambda\times\Lambda)=L^2(\Lambda)\ot L^2(\Lambda)$ we define
an operator $R:f(x,y)\ra e^{ir(x,y)}f(y,x)$. It follows easily that 
$R\in\cb(L^2(\Lambda\times\Lambda))$ is such that $R^2=\idt$, $R=R^+$ and
$R$ obeys the YB condition. In the space $L^2(\Lambda^{\ot n})$ we define a 
representation of the permutation group $S_n$ be defining the action of
the transposition 
$\tau_i:f(x_1,...,x_n)=e^{ir(x_i,x_{i+1})}f(x_1,...,x_{i+1},x_i,...,x_n)$ and 
then we extend this definition to arbitrary permutation by taking the product
of the corresponding transpositions. In this way we obtain a homomorphism 
$S_n\ni\pi\ra R(\pi)\in\cb(L^2(\Lambda^{\times n}))$ in such a way that the 
operator $R_n=\frac{1}{n!}\sum_{\pi\in S_n}R_n(\pi)$ is an orthogonal 
projection.
Let $\Gamma(L^2(\Lambda))$ be the free Fock space over $L^2(\Lambda)$ and
let $R=\bigoplus_{n=0}^\infty R_n$ be the corresponding projector. By a $r$-Fock
module we mean the space $P_R\Gamma(L^2(\Lambda))$. Let $b(f)$ $(b^*(f))$
be the standard anihilation (resp. creation) operators acting the free space
$\Gamma(L^2(\Lambda))$ and let us define $a_r(f)=P_Rb(f)P_R\equiv\int dx 
f(x)a(x)$ 
and $a_r^+(f)=P_Rb^*(f)P_R\equiv\int dx f(x)a^+(x)$. By simple calculations 
we get:
$$
a_r(x)a_r^+(y)-e^{ir(x,y)}a_r^+(y)a_r(x)=\delta(x-y)
$$
$$
a_r(x)a_r(y)-e^{ir(x,y)}a_r(y)a_r(x)=0
$$
$$
a_r^+(x)a_r^+(y)-e^{ir(x,y)}a_r^+(y)a_r^+(x)=0
$$
on a suitable domain in $\Gamma_r(L^2(\Lambda))$.
The main observation is the following:
\begin{twierdz}
Let $(\Lambda_n,\sigma_n)$ be a family of bounded regions as above and 
such that $\Lambda_n\subset \Lambda_{n+1}$ and $\bigcup_n\Lambda_n=\Bbb{R}^d$
and a family $\sigma_n\in C^2(\partial\Lambda_n)$, $\sigma_n\ge0$ for all n.
Assume that $r(x-y)$ is such that $\nabla r\in L^2(\Bbb{R}^d)$ and
$\nabla r(x)\nabla r(x-z)\in L^1(\Bbb{R}^{2d})$.

Define the finite volume free energy density
$$
P_r^{(\Lambda_n,\sigma_n)}(\beta,\mu)=
\frac{1}{|\Lambda_n|}ln\;Tr_{\Gamma_r(\Lambda_n)}(e^{-\beta h^{n}})
$$
where $h^n=-\Delta_{\Lambda_n}^{\sigma_n}-\mu \idt$.
Then there exists $(z_0,\beta_0)$ such that for all $z\in\Bbb{C}:|z|<z_0$,
$\beta<\beta_0$ there exists
$$
\lim_{n\ra\infty}P_r^{(\Lambda_n,\sigma_n)}(\beta,\mu)\equiv P_r(\beta,\mu)
$$
and the limiting free energy density does not depends on $(\sigma_n)$ and 
$(\Lambda_n)$.
Moreover $P_r$ is an analitic in the circle $|e^{-\beta\mu}\equiv z|<z_0$.
\end{twierdz}

Outline of the proof:

{\bf Step 1} Let $T_+(x_1,x_2)\equiv e^{-\frac{i}{2}r(x_1,x_2)}$. Then the 
map 
$
U_+(R)f_n(x_1,...,x_i,...,x_j,...,x_n)
\equiv\prod_{i,j=1}^nT_+(x_i,x_j)f_n(x_1,..,x_j,...,x_i,...,x_n)
$
extended to $\Gamma_R(\cch)$ gives rise to an unitary isomorphism
between $\Gamma_R(\cch)$ and $\Gamma_+(\cch)$, so that (see also~\cite{LM}):
$$
H_+^{(\Lambda,\sigma)}\equiv U_+(R)d\Gamma_R(-\Delta^{(\Lambda,\sigma)}-\mu)
U_+^{-1}(R)=d\Gamma_+(-\Delta^{(\Lambda,\sigma)}-\mu)+R_+^{(\Lambda,\sigma)}
$$
where
\beqna*
R_+^{(\Lambda,\sigma)}
&=&\int_\Lambda a^+(x)a^+(y)\{ \frac{1}{4}(\nabla_xr(x-y))^2\} dxdy\\
&&+\int_\Lambda a^+(x)a^+(y)a^+(z)
\{ \frac{1}{12}\nabla_xr(x-y)\nabla_yr(y-z)\\
&&+\mbox{ cyclic permutations }\}a(x)a(y)a(z)dxdydz
\eeqna*
{\bf Step 2} Let $dW_{x,y}^{(\Lambda,\sigma),\beta}$ be the corresponding 
Wiener bridge measure on the space $C([0,\beta]\ra\overline{\Lambda})$.
Then be using step 1 we have:

\beqna*
Tr_{\Gamma_T(\cch)}(e^{-\beta(\Delta^{(\Lambda,\sigma)}-\mu)})
&=&Tr_{\Gamma_+(\cch)}e^{-\beta H_+^{(\Lambda,\sigma)}}\\
&=&\sum_{n=0}^\infty\frac{1}{n!}\sum_{j_1=1}^\infty\frac{1}{j_1}...
\sum_{j_n=1}^\infty\frac{1}{j_n}
\int_\Lambda dx_1\int dW_{x_1|x_1}^{(\Lambda,\sigma),j_1\cdot\beta}(w_1)\\
&&...
\int_\Lambda dx_n\int dW_{x_n|x_n}^{(\Lambda,\sigma),j_n\cdot\beta}(w_n)
exp-\beta \ce^\beta(w_1,...,w_n)
\eeqna*
where $\ce^\beta$ is the corresponding energy factor (see i.e.~\cite{BR},
p 381 and also~\cite{GR}).

In this step we generalise the analysis of Ginibre~\cite{Ginibre} to the case
multibody potentials, (see~\cite{GR} for details). In particular we obtain
the existence and analicity of
$\lim_{n\ra\infty}P_r^{(\Lambda_n,\sigma_n)}(\beta,\mu)=P_r(\beta,\mu)$
for $\beta$ and $\mu$ as stated.

{\bf Step 3} In this step we prove the indenpedence of 
$P^{(\Lambda_n,\sigma_n)}(\beta,\sigma)$ on $(\Lambda_n,\sigma_n)$ by adopting 
the
methods of~\cite{GR}. 

\begin{uwaga}
For a complete proof we refer to~\cite{GR} where also some generalisations and
also complementary results are included.

Let $d=2$ and let $u$ be fixed unit vector in the plane $\Bbb{R}^2$.
Taking the formfactor $r$ of the following form:
$$
r(x,y)=-\theta sign[(x-y)_k\cdot \t{u}_k]
$$
where $\theta\in\Bbb{R}$ is the so called statistical parameter and $\t{u}$ is
the dual of $u$,
it follows that the corresponding hamiltonian $H_+^{(\Lambda,\sigma)}$ (properly 
defined)
coincidies with  that of Leinaas-Myrheim~\cite{Leinaas}.
\end{uwaga}

\section{Thermal states on $q$-commutation relations algebras}
Let $\cw(\cch ,T)$ be a Wick algebra with positive Fock 
module $\Gamma_T(\cch)$. Let $U_t=e^{ith}$ be one parameter weakly 
continous unitary group on $\cch$. The second quantisation $\Gamma_T(U_t)$
of $U_t$ forms unitary group on $\Gamma_T(\cch)$ iff $[d\Gamma_2(h),T]=0$.
In particular for q relations any unitary group  $U_t$ on $\cch$ lifts to 
an unitary group $\Gamma_q(U_t)$ on $\Gamma_q(\cch)$ which moreover gives
one parameter weakly continous group $\alpha_t$ of $*$-automorphisms of the
${\bf W}^*$-algebra $\mbox{\gt M}_q$ defined as weak closure of the $*$-algebra 
generated by the
Fock realisation of the corresponding $q$-commutation relation operators.
Let us recall that a linear continous functional $\omega_\beta$ on $\mbox{\gt 
M}_q$ 
is called $\beta$-KMS functional
for a given $\alpha_t$ iff any $A,B\in \mbox{\gt M}_q$ there exists a function 
$F_{A,B}(z)$
holomorphic in the strip $0<Imz<\beta$ and continous on the boundaries 
$\{ z\in \Bbb{C}:\;\;Imz=0\}\cup\{ z\in\Bbb{C}:\;\;Imz=\beta\}$ and such that:
$$
F_{A,B}(t)=\omega(A\alpha_t(B));\;\;\;\;F_{A,B}(t+i\beta)=\omega(\alpha_t(B)A)
$$
It is interesting to note that althought the thermodynamics of the quon matter 
(for $q\in(-1,1)$) does not depends on $q$, the corresponding KMS states  
$\omega_\beta$
explicitly depend on q.
\begin{stwierdz}\label{p1}
Let $q\in(-1,1)$ and let $\mbox{\gt M}_q$ be the corresponding (Fock 
realisation) 
of ${\bf W}^*$ algebra of $q$-commutation relation. Assume that $h=h^+$ is 
strictly
positive selfadjoint generator of $U_t=e^{ith}$. If $\omega_\beta$ is
$\beta$-KMS state on $\mbox{\gt M}_q$ then
\be\label{l1}
\omega_\beta^{(2)}(a(f)a^+(g))=\<f,\frac{1}{1-qe^{-\beta h}}g\>
\ee
\be\label{l2}
\omega_\beta^{(2)}(a^+(f)a(g))=\< g,\frac{e^{-\beta h}}{1-qe^{-\beta h}}f\>
\ee
for any $f,g\in\cch$
\end{stwierdz}
\begin{Proof}
By passing to the analytic elements of $\alpha_t$ and the use 
of the $q$-commutation relations, see also~\cite{GR}.
\end{Proof}

The next problem is to extend the formulas \ref{l1} and \ref{l2} to the whole
algebra $\mbox{\gt M}_q$. Keeping in mind the form of the vacuum functional for 
the Fock
representation of the $q$-commutation relations~\cite{BS1,BS2} we can deduce the 
validity 
of the following proposition:
\begin{stwierdz}
Let $q\in(-1,1)$, $\mbox{\gt M}_q, U_t=e^{ith}$ be as in Proposition \ref{p1}. 
Then the
functional $\omega_\beta$ defined by:
\be\label{l3}
\omega_\beta(a_{f_1}^{\sigma_1},..,a_{f_n}^{\sigma_n})=
\left\{
\ba{ccc}
0&\mbox{for}&n=2k+1\\
\sum_{\mu\in P_2(k)
\subset\prod^c(n),\alpha_i<\beta_i
}
q^{\#\mu}\prod_{j=1}^k\omega_\beta^{(2)}(a_{f_{\alpha_j}}^{\sigma_{\alpha_j}}
a_{f_{\beta_j}}^{\sigma_{\beta_j}})&\mbox{for}&n=2k
\ea
\right.
\ee
where 
$$
P_2(k)=\{((\alpha_1,\beta_1),...,(\alpha_k,\beta_k))|\;\alpha_i<\beta_i,\;\;
i=1,...,k,\;a_f^{-1}\equiv a_f,\; a_f^{+1}\equiv a_f^+\}
$$ 
and $\#\mu$ is the number of crossing in two-partitions $\mu\in P_2(k)$ and 
$\omega_\beta^{(2)}$
is given by \ref{l1}, \ref{l2}, gives $\beta$-KMS functional on $\mbox{\gt M}_q$ 
with
respect to $\alpha_t$ for any $q\in[-1,1]$.
\end{stwierdz}
The functional $\omega_\beta$ given by \ref{l3} is gauge-invariant and 
has the property that 
$\lim_{\beta\ra\infty}\omega_\beta=\mbox{ Fock vacuum functional}$.
We conjecture:
\begin{hipoteza}
The functional $\omega_\beta$ given by \ref{l3} and additionally 
$\omega_\beta(\idt)=1$ is a state on $\mbox{\gt M}_q$.
\end{hipoteza}
\begin{uwaga}
The interesting feature of the above conjecture is that if it is true,
then $\omega_\beta$ leads to a class of the new representations 
of the $q$-commutation relations. For $q=0$ the validity of this conjecture 
has been proven in~\cite{Shl} and for $q\in\{-1,1\}$ the formula \ref{l3} gives 
a general form of quasi-free gauge-invariant $\beta$-KMS functional on the
corresponding CAR and CCR algebras.
\end{uwaga}


\end{document}